\documentclass{SCIYA2017enOL}
\online
\begin{document}

\ensubject{fdsfd}%二级学科

%%%%%%%%%%%%%%%%%%%%%%%%%%%%%%%%%%%%%%%%%%%%%%%%%%%%%%%
%%% Authors do not modify the information below
%%% 作者不需要修改此处信息
%%% 有专题名称时, 将第一行的{}注释掉, 使用第二行
\ArticleType{ARTICLES}%栏目
%\SpecialTopic{Progress of Projects Supported by NSFC}%专题
%\SubTitle{Dedicated to Professor Yang Lo on the Occasion of his {\rm 70}th Birthday}%专刊说明
\Year{2018}
\Month{December}%
\Vol{60}
\No{1}
\BeginPage{1} %
\DOI{10.1007/s11425-000-0000-0}
\ReceiveDate{March 27, 2018}
\AcceptDate{December 10, 2018}
%\OnlineDate{January 1, 2017}
%%%%%%%%%%%%%%%%%%%%%%%%%%%%%%%%%%%%%%%%%%%%%%%%%%%%%%%

%%% title: 标题
%%%   \title{title}{title for citation}
\title[]{Nonexistence of positive supersolution  to a class of semilinear
elliptic equations and systems in an exterior domain}
{Nonexistence of positive supersolution  to a class of semilinear
elliptic equations and systems in an exterior domain}%%后边花括号是眉题

%%% Corresponding author: 通信作者
%%%   \author[number]{Full name}{{email@xxx.com}}
%%% General author: 一般作者
%%%   \author[number]{Full name}{}
\author[1]{Huyuan Chen}{chenhuyuan@yeah.net}
\author[1,$\ast$]{Rui Peng}{pengrui$\_$seu@163.com}
\author[2]{Feng Zhou}{fzhou@math.ecnu.edu.cn}

%%% Author information for page head. 页眉中的作者信息
%%% 若此处指定以此处为准, 否则直接调用author信息
\AuthorMark{Huyuan Chen}

%%% Authors for citation. 首页引用中的作者信息
%%% 若此处指定以此处为准, 否则直接调用author信息
\AuthorCitation{Huyuan Chen, Rui Peng, Feng Zhou}

%%% Address. 地址
%%%   \address[number]{Address, City {\rm Postcode}, Country}
\address[1]{School of Mathematics and Statistics, Jiangsu Normal University,\\
Xuzhou, Jiangsu {\rm 221006}, China}
\address[2]{Center for PDEs and Department of Mathematics, East China Normal University,\\
 Shanghai {\rm 200241}, China}

%%% Abstract. 摘要
\abstract{In this paper, we primarily consider the following semilinear elliptic equation
 \begin{eqnarray*}
 \arraycolsep=1pt\left\{
 \begin{array}{lll}
 \displaystyle  -\Delta u= h(x,u)\quad \
 &{\rm in}\ \Omega,\\[1.5mm]
 \phantom{   -\Delta }
 \displaystyle  u\ge 0\qquad  &{\rm   on}\ \partial{\Omega},
 \end{array}\right.
 \end{eqnarray*}
where $\Omega$ is an exterior domain in $\mathbb{R}^N$ with $N\ge 3$,   $h: \Omega\times \mathbb{R}^+\to \mathbb{R}$ is a measurable function,  and
derive optimal nonexistence results of positive supersolution. Our argument is
based on a nonexistence result of positive supersolution of a linear elliptic
problem with Hardy potential. We also establish sharp nonexistence results of positive supersolution
to an elliptic system.}

%%% Keywords. 关键词
\keywords{Semilinear elliptic problem; Supersolution; Nonexistence}

\MSC{35J60, 35B53}

\maketitle

%%%%%%%%%%%%%%%%%%%%%%%%%%%%%%%%%%%%%%%%%%%%%%%%%%%%%%%
%%% The main text. 正文部分%
%%  图表引用\cref公式引用\eqref参考文献\cite
%%%%%%%%%%%%%%%%%%%%%%%%%%%%%%%%%%%%%%%%%%%%%%%%%%%%%%%
\section{Introduction}

In this paper, we are mainly concerned with the nonexistence of positive
supersolution to the following semilinear elliptic equation
 \begin{equation}\label{eq 1.1}
 \arraycolsep=1pt\left\{
 \begin{array}{lll}
 \displaystyle  -\Delta u= h(x,u)\quad \
 &{\rm in}\ \Omega,\\[1.5mm]
 \phantom{   -\Delta }
 \displaystyle  u\ge 0\qquad  &{\rm   on}\ \partial{\Omega},
 \end{array}\right.
\end{equation}
where $\Omega$ is a punctured or an exterior domain in $\mathbb{R}^N$ with $N\ge3$;
that is, $\Omega=\mathbb{R}^N\setminus \mathcal{O}$, here $\mathcal{O}$ is
a bounded, closed smooth subset of $\mathbb{R}^N$ and $h: \Omega\times \mathbb{R}^+\to \mathbb{R}$ is a measurable function, here $\mathbb{R}^+=[0,+\infty)$.
We assume, unless otherwise specified, that $N\ge3$ throughout the paper.
Without loss of generality, we also assume that $\mathcal{O}\subset \overline B_R(0)$,
where $\overline B_R(0)$ represents the ball with the radius $R$, centered at the origin.
A typical punctured domain is $\Omega=\mathbb{R}^N\setminus \{0\}$, and a typical exterior domain is
$\Omega=\mathbb{R}^N\setminus \overline B_R(0)$ with $\mathcal{O}=\overline B_R(0)$.
{\it A function $u\in C^{2}(\Omega)$ is said to be a positive supersolution of (\ref{eq 1.1}) if $u(x)>0$ and
$-\Delta u(x)\ge h(x,u)$ for all $x\in\Omega$.  }

The existence and nonexistence of solution or supersolution to problem (\ref{eq 1.1}) have attracted great attentions for many years; see \cite{AMQ,CA,A,AS,BK,C,CF,DLY,DG,GY,L1,L2,LLd} and the references therein. In the special case that $h$ only depends on
$u$, given $R_0>0$, Alarcon, Melian and Quaas \cite{AMQ} proved that problem (\ref{eq 1.1}) with $\Omega=\mathbb{R}^N\setminus B_{R_0}(0)$ admits a positive solution {\it if and only if} $h$ satisfies
 $$
 \int_0^{\sigma_0}h(t)t^{-\frac{2N-2}{N-2}}dt<+\infty
 $$
for some $\sigma_0>0$, by appealing to ODE techniques. Nevertheless, such an approach fails to apply to (\ref{eq 1.1})
if the nonlinear term $h$ is not radially symmetric.

In particular, when $h(x,u)=V(x)u^p$, (\ref{eq 1.1}) becomes the following
 \begin{equation}\label{eq 1.1-aa}
 \arraycolsep=1pt\left\{
 \begin{array}{lll}
 \displaystyle  -\Delta u=V(x)u^p\quad \
 &{\rm in}\ \Omega,\\[1.5mm]
 \phantom{   -\Delta }
 \displaystyle  u\ge 0\qquad  &{\rm   on}\ \partial{\Omega}.
 \end{array}\right.
 \end{equation}
The traditional method to establish
nonexistence results for solution or supersolution to (\ref{eq 1.1-aa}) is to
make use of the fundamental solution and Hadamard property (\cite{AS,B,BP,CF}).
It is worth mentioning that if $p =-2$,  problem (\ref{eq 1.1-aa}) in a bounded domain
is used to describe the MEMS model (\cite{EGG,G}),
and if $p =-1$, problem (\ref{eq 1.1-aa}) is related to the study of singular minimal
hypersurfaces with symmetry (\cite{M,S}).

When $V(x)=(1+|x|)^{\beta}$ and $\Omega$ is a punctured or an exterior domain, Bidaut-V\'{e}ron \cite{B}
and Bidaut-V\'{e}ron and Pohozaev \cite{BP} showed that
problem (\ref{eq 1.1-aa}) has no solution
if
 $$
 p\le \frac{N+\beta}{N-2}:=p_\beta^*\ \ \ \ \mbox{with}\ \beta\in(-2,2).
 $$
On the other hand, when $\Omega=\mathbb{R}^N\setminus\{0\}$ and
$V(x)=|x|^{a_0}(1+|x|)^{\beta-a_0}$ with $a_0\in(-N,+\infty)$ and $\beta\in(-\infty,a_0)$,
Chen, Felmer and Yang \cite{CFY} derived infinitely many positive solutions of (\ref{eq 1.1-aa})
if
 $$
 p\in\left(p_\beta^*,\frac{N+a_0}{N-2}\right)\cap (0,\,+\infty).
 $$

Armstrong and Sirakov \cite{AS} and Chen and Felmer \cite{CF} dealt with the more general potential
 $$
 V(x)\ge |x|^{\beta} (\ln |x|)^\tau  \quad{\rm for}\  |x|>e,
 $$
where $\beta>-2$ and $\tau\in  \mathbb{R}$. Especially, \cite[Theorem 3.1]{AS} and \cite[Theorem 1.1]{CF}
imply the following  result:

 \begin{itemize}
 \item[$\bullet$] {\it Let  $\beta\in(-2,2)$. Problem (\ref{eq 1.1-aa}) with $\Omega=\mathbb{R}^N\setminus B_e(0)$
 has no positive supersolution provided that either
 $p < p_\beta^*,\,\tau\in \mathbb{R}$ or $p=p_\beta^*,\,\tau\ge0$.}
 \end{itemize}

In the current paper, we will provide a sharp improvement of the above result;
indeed, we can conclude the following:

\begin{itemize}
 \item[$\bullet$] {\it Let  $\beta>-2$.  Problem (\ref{eq 1.1-aa}) with $\Omega=\mathbb{R}^N\setminus B_e(0)$
 has no positive supersolution provided that either
 $1\leq p < p_\beta^*,\,\tau\in\mathbb{R}$ or $p=p_\beta^*,\,\tau>-1$ (see Proposition \ref{cr 2});}

\item[$\bullet$] {\it Let  $\beta=-2$ and  $p=p_\beta^*=1$. Problem (\ref{eq 1.1-aa}) with $\Omega=\mathbb{R}^N\setminus B_e(0)$
 has no positive supersolution provided that
  $\liminf_{|x|\to+\infty} V(x)|x|^2>\frac{(N-2)^2}{2}$ (see Theorem \ref{teo 2}); }

 \item[$\bullet$] {\it Let  $\beta\in\mathbb{R}$. Problem (\ref{eq 1.1-aa}) with $\Omega=\mathbb{R}^N\setminus B_\ell(0)$ has positive supersolution for properly large $\ell$ provided that either
 $p>p_\beta^*,\,\tau\in\mathbb{R}$ or $p=p_\beta^*,\,\tau<-1$ (see Proposition \ref{cr 2-5}).}

 \end{itemize}

As a consequence, the above results show that both $p=p_\beta^*$ and $\tau=-1$ are the critical values for the existence
of positive supersolution to (\ref{eq 1.1-aa}).

For the more general nonlinear problem (\ref{eq 1.1}), let us assume that
 $$
 h(x,u)\geq\tilde h(x,u)\geq0\ \ \mbox{for all}\ x\in\mathbb{R}^N\setminus B_e(0),\,u\geq0,
 $$
and $\tilde h:\Omega\times \mathbb{R}^+\to\mathbb{R}^+$ is a function  satisfying  the following
 \begin{itemize}
\item[$\rm(H)$]
(a) for any $x\in \mathbb{R}^N\setminus B_e(0)$,
 $$
 \frac{\tilde h  (x,s_1)}{s_1}\ge \frac{\tilde h  (x,s_2)}{s_2}\quad{\rm if}\  s_1\ge s_2>0;
 $$

(b) for any $t>0$,
 $$
 \liminf_{|x|\to+\infty} \tilde h  (x,\,t|x|^{2-N})  |x|^{N}>\frac{(N-2)^2}{4};
 $$

if (b) fails, we assume that

 (b1) there exists  $\sigma_1\in(0,1)$ such that for any $t>0$,
 $$
 \liminf_{|x|\to+\infty} \tilde h  (x,\,t|x|^{2-N})  |x|^{N}  (\ln|x|)^{\sigma_1}>0;
 $$
and

(b2) there exists  $\sigma_2>0$ such that for any $t>0$,
 $$
 \liminf_{|x|\to+\infty}\frac{\tilde h  (x, t |x|^{2-N} (\ln|x|)^{\sigma_2}) }{t |x|^{-N} (\ln|x|)^{\sigma_2}  }>\frac{(N-2)^2}{4}.
 $$
\end{itemize}
Then we have

\begin{theorem}\label{teo 1} Under the assumption {\rm (H)}, problem (\ref{eq 1.1})  has no positive supersolution.

\end{theorem}

We would like to mention that (H)-(a) means that $h$ is linear or superlinear,
(H)-(b) is related to the subcritical case while (H)-(b1)(b2) deal with the critical case. As one will see below,
Theorem \ref{teo 1} allows us to obtain some optimal nonexistence results of positive supersolution to problem (\ref{eq 1.1-aa}).

To prove Theorem \ref{teo 1}, it turns out that a nonexistence result of positive supersolution of
the linear Hardy elliptic problem (\ref{eq 1.1-ab})
(see Section 2) is vital in our analysis. Such a nonexistence result can be established by Agmon-Allegretto-Piepenbrink theory \cite{As}; in this paper, we shall provide a different proof which seems simpler.

\smallskip Another focus of our paper is on the following semilinear elliptic system:
 \begin{equation}\label{sys 1.1}
 \arraycolsep=1pt\left\{
 \begin{array}{lll}
 \displaystyle  -\Delta u= h_1(x,u,v)\qquad
 &{\rm in}\  \Omega,\\[1.5mm]
 \displaystyle  -\Delta v= h_2(x,u,v)\qquad
 &{\rm in}\  \Omega,\\[1.5mm]
 \phantom{   -\Delta }
 \displaystyle  u,v\ge 0\qquad  &{\rm   on}\ \partial{\Omega}.
 \end{array}\right.
 \end{equation}
It is known that Liouville-type theorems for system (\ref{sys 1.1}) have been established (see \cite{DF,FG,QS,SZ,Souto}) primarily on the whole space $\Omega=\mathbb{R}^N$. In particular, when $h_1(x,u,v)=v^p,\ h_2(x,u,v)=u^q$, the nonexistence of positive solution to (\ref{sys 1.1}) in $\Omega=\mathbb{R}^N$ has been investigated  by \cite{DF,SZ,Souto}
in the subcritical case that
 $$
 \frac{1}{p+1}+\frac{1}{q+1}>1-\frac2N.
 $$

It seems that there is little research work devoted to the nonexistence/existence of solution/supersolution to (\ref{sys 1.1}) when $\Omega$ is an exterior domain. {\it We say a function pair $(u,v)\in C^{2}(\Omega)\times C^{2}(\Omega)$ is a positive supersolution of (\ref{sys 1.1}) if $u(x),\,v(x)>0$, $-\Delta u(x)\ge h_1(x,u,v)$ and
$-\Delta v\ge h_2(x,u,v)$ for all $x\in\Omega$.}

\vskip5pt
For system (\ref{sys 1.1}),  two functions $\tilde h_1$, $\tilde h_2$ are involved to control functions $h_1,\ h_2$ respectively.

 \begin{itemize}
\item[$\rm (SH)$] Suppose that the nonnegative function $\tilde h_i\, (i=1,2)$ defined on $(\mathbb{R}^N\setminus B_e(0))\times[0,+\infty)\times[0,+\infty)$ satisfies that

\vskip5pt
(a) the maps $t\to \tilde h_1(x,s,t)$, $s\to \tilde h_2(x,s,t)$ are nondecreasing and
for any $x\in\mathbb{R}^N\setminus B_e(0)$,
 $$
 \frac{\tilde h_1(x,s_1,t_1)}{s_1}\ge \frac{\tilde h_1(x,s_2,t_2)}{s_2}\quad{\rm and}\quad\frac{\tilde h_2(x,s_1,t_1)}{t_1}\ge
\frac{\tilde h_2(x,s_2,t_2)}{t_2}\quad{\rm if}\  s_1\ge s_2>0,\ t_1\ge t_2>0;
$$

(b) for any $t>0$,
 $$
 \liminf_{|x|\to+\infty} \tilde h_1(x,\,t|x|^{2-N},\, t|x|^{2-N})  |x|^{N}>\frac{(N-2)^2}{4},
 $$
or
 $$
 \liminf_{|x|\to+\infty} \tilde h_2(x,\,t|x|^{2-N},\, t|x|^{2-N})  |x|^{N}>\frac{(N-2)^2}{4};
 $$

if (b) fails, we assume that

 (b1) there exist positive constants $\sigma_3, \sigma_4$ with either $\sigma_3<1$ or $\sigma_4<1$,
 such that for any $ t>0$,
 $$
 \liminf_{|x|\to+\infty} \tilde h_i(x,\,t|x|^{2-N},\, t|x|^{2-N})  |x|^{N}   (\ln |x|)^{\sigma_{i+2}}>0,\quad i=1,2;
 $$
and

 (b2) if $\sigma_3<1$, $\sigma_4\ge 1$
 there exists $ \sigma_5>0$  such that for any $t>0$,
 $$
 \liminf_{|x|\to+\infty}\frac{\tilde h_1(x,\, t|x|^{2-N}
 (\ln |x|)^{\sigma_5}, t|x|^{2-N} ) }{t|x|^{-N} (\ln |x|)^{\sigma_5}}>\frac{(N-2)^2}{4};
 $$
  if $\sigma_3\ge 1$, $\sigma_4<1$, there exists $ \sigma_6>0$  such that for any $t>0$,
 $$
 \liminf_{|x|\to+\infty}\frac{\tilde h_2(x,\, t|x|^{2-N},t|x|^{2-N}
 (\ln |x|)^{\sigma_6})) }{t|x|^{-N} (\ln |x|)^{\sigma_6}  }>\frac{(N-2)^2}{4};
 $$
 (b3) if $\sigma_3,\, \sigma_4<1$, there exist $\sigma_5, \sigma_6>0$  such that
 for any $t>0$,
 $$
 \liminf_{|x|\to+\infty}\frac{\tilde h_1(x,\, t|x|^{2-N} (\ln |x|)^{\sigma_5}, t|x|^{2-N}
 (\ln |x|)^{\sigma_6}) }{t|x|^{-N} (\ln |x|)^{\sigma_5} }>\frac{(N-2)^2}{4}
 $$
 or
 $$
 \liminf_{|x|\to+\infty}\frac{\tilde h_2(x,\, t|x|^{2-N} (\ln |x|)^{\sigma_5},t|x|^{2-N}
 (\ln |x|)^{\sigma_6})) }{t|x|^{-N} (\ln |x|)^{\sigma_6}  }>\frac{(N-2)^2}{4}.
 $$

\end{itemize}
Then we can state

\begin{theorem}\label{teo 3} Assume that for all $(x,u,v)\in (\mathbb{R}^N\setminus B_e(0))\times [0,+\infty)\times [0,+\infty)$,
the functions $h_1,\,h_2$ satisfy
 $$
 h_1(x,u,v)\ge \tilde h_1 (x,u,v), \ \quad h_2(x,u,v)\ge \tilde h_2 (x,u,v),
 $$
with $\tilde h_1, \tilde h_2$ fulfilling {\rm(SH)}. Then system (\ref{sys 1.1})  has no positive supersolution.

\end{theorem}

Assumption (SH)-(b) is related to the subcritical case,  (SH)-(b2) says that one of the nonlinearities is critical, and (SH)-(b3) represents that both nonlinear terms are critical. In particular, when the nonlinearities $h_1,\ h_2$ take the form  $h_1(x,u,v)=|x|^{\beta_1} (\ln |x|)^{\tau_1}u^{p_1}v^{q_1}$ and $h_2(x,u,v)=|x|^{\beta_2} (\ln |x|)^{\tau_2}u^{p_2}v^{q_2}$, we  are able to clarify the nonexistence and existence of positive supersolution in terms of the parameters $p_1,p_2,q_1,q_2,\tau_1$ and $\tau_2$; see Proposition \ref{cr 2-3} and Proposition \ref{cr 2-6-}  for the precise details.

The rest of the paper is organized as follows. In Section 2, we show the nonexistence of supersolution of a linear Hardy problem. In Section 3, we prove our main results  Theorems \ref{teo 1} and \ref{teo 3}. In Section 4, we apply Theorems \ref{teo 1} and \ref{teo 3} to
two concrete examples to obtain sharp nonexistence results.

\setcounter{equation}{0}
\section{Nonexistence of positive supersolution of a linear Hardy problem}

In this section,  we shall investigate the linear elliptic problem with Hardy potential:
 \begin{equation}\label{eq 1.1-ab}
 \arraycolsep=1pt\left\{
 \begin{array}{lll}
 \displaystyle  -\Delta u=V(x)u\quad \
 &{\rm in}\ \Omega,\\[1.5mm]
 \phantom{   -\Delta }
 \displaystyle  u\ge 0\qquad  &{\rm   on}\ \partial{\Omega}.
 \end{array}\right.
 \end{equation}

The nonexistence result of positive supersolution to (\ref{eq 1.1-ab}) reads as follows.

 \begin{theorem}\label{teo 2} Assume that $\Omega$ is a punctured or an exterior domain, and
 $V$ is a nonnegative function satisfying
 \begin{equation}\label{v2}
 \liminf_{|x|\to+\infty} V(x)|x|^2>\frac{(N-2)^2}{4}.
 \end{equation}
Then problem (\ref{eq 1.1-ab}) has no positive supersolution.

\end{theorem}

We remark that ${(N-2)^2}/{4}$ is the best constant in the Hardy-Sobolev inequality
 $$
 \int_{\mathbb{R}^N}|\nabla u|^2 dx \ge \frac{(N-2)^2}{4}\int_{\mathbb{R}^N}\frac{u^2}{|x|^2}dx.
 $$
Related Hardy problems have been studied extensively; one may refer to, for instance, \cite{BDT,BV1,CC,CQZ,D}. In particular,
the authors \cite{CQZ} considered the Hardy problem
$-\Delta u=\frac\mu{|x|^2} u+f(x)$ in $\mathcal{D}\setminus\{0\}$, subject to the zero Dirchlet boundary condition,
where $f\geq0$ and $\mathcal{D}$ is a bounded domain containing the origin, and proved that this problem has no positive solution once $\mu>\frac{(N-2)^2}{4}$. To this end, they used new distributional identities
to classify the isolated singular solution of $-\Delta u=\frac\mu{|x|^2} u+f$ in  $\mathcal{D}\setminus\{0\}$
and find fundamental solutions of  $-\Delta u=\frac\mu{|x|^2} u$ in $\mathbb{R}^N$.

It is worth noting that \cite{DFP,P,P1} indicate that the nonexistence of positive supersolution to (\ref{eq 1.1-ab}) can be obtained by using Agmon-Allegretto-Piepenbrink theory \cite{As}. We will provide a different and elementary proof. Our strategy is to employ the Kelvin transform to transfer the unbounded domain $\Omega$ into a bounded one containing the origin.

For the linear elliptic equation  involving the general homogeneous potential in the punctured domain $\mathbb{R}^N\setminus\{0\}$:
 \begin{equation}\label{eq 1.1-p1}
 \displaystyle  -\Delta u=\mu |x|^{-\alpha} u\quad \
 {\rm in}\ \, \mathbb{R}^N\setminus\{0\},
 \end{equation}
 we obtain that

\begin{theorem}\label{coro 1} Problem (\ref{eq 1.1-p1}) has no positive supersolution provided that one of the following conditions holds:

 \begin{itemize}
 \item[{\rm (i)}] $\alpha\not=2,\ \mu>0;$

 \item[{\rm (ii)}] $\alpha=2,\ \mu>\frac{(N-2)^2}{4}.$

 \end{itemize}

\end{theorem}

Theorem \ref{coro 1} is optimal in a certain sense, and it also reveals essential differences between problem (\ref{eq 1.1}) with a
punctured domain and problem (\ref{eq 1.1}) with an exterior domain; see the following remark.

 \begin{remark}\label{re-a} Concerning Theorem \ref{coro 1}, we would like to make some comments as follows.

 \begin{itemize}
 \item[{\rm (i)}] When $\alpha<2,\ \mu>0$, one can easily see from the proof of Theorem \ref{coro 1} that
 the linear problem
 \begin{equation}\label{eq 1.1-p2}
 \arraycolsep=1pt\left\{
 \begin{array}{lll}
 \displaystyle  -\Delta u=\mu |x|^{-\alpha} u\quad \
 &{\rm in}\ \, \mathbb{R}^N\setminus B_\ell(0),\\[1.5mm]
 \phantom{   -\Delta }
 \displaystyle  u\ge 0\qquad  &{\rm   on}\ \, \partial{B_\ell(0)}
 \end{array}\right.
 \end{equation}
 has no positive supersolution for any $\ell>0;$

 \item[{\rm (ii)}] When $\alpha>2,\ \mu>0$, problem (\ref{eq 1.1-p2}) has positive supersolution for properly large $\ell$; this can be seen from Proposition \ref{cr 2-5}(i) below (by taking $p=1,\,\beta<-2,\,\tau=0$ there). Such a result is in sharp contrast with the above (i) and Theorem \ref{coro 1}(i)$;$

  \item[{\rm (iii)}]  When $\alpha=2,\ \mu\leq {(N-2)^2}/{4}$, problem (\ref{eq 1.1-p1}) (and so problem (\ref{eq 1.1-p2})) has
positive supersolution by considering the fundamental solutions of Hardy operators; one can refer to, for example, \cite{CQZ}. Hence, for $\alpha=2$, $\mu={(N-2)^2}/{4}$ is the critical value of existence of positive supersolution of problems (\ref{eq 1.1-p1}) and (\ref{eq 1.1-p2}).

 \end{itemize}

\end{remark}

\subsection{Proof of Theorem \ref{teo 2}}

In the case that $\beta=-2$ and so the critical exponent $p^*_\beta=1$, problem (\ref{eq 1.1-ab})
is related to the Hardy-Leray potentials.  The following result plays an essential  role
in the proof of Theorem \ref{teo 2}.

\begin{lemma}  \cite[Proposition 5.2]{CQZ}  \label{sub}
  Assume that $\mu>(N-2)^2/4$,  $\mathcal{D}$ is a bounded smooth domain containing the origin and
   $f\in L^\infty_{loc}(\overline{\mathcal{D}\setminus \{0\}})$ is a nonnegative function. Then the Hardy problem
\begin{equation}\label{eq 2.2}
 \arraycolsep=1pt\left\{
\begin{array}{lll}
 \displaystyle -\Delta u=\frac{\mu}{|x|^2}u+f\quad
   &{\rm in}\  \mathcal{D}\setminus\{0\} ,\\[1.5mm]
 \phantom{   -\Delta }
 \displaystyle  u=0 \qquad  &{\rm   on}\   \partial\mathcal{D}
 \end{array}\right.
 \end{equation}
 has no positive solution.
\end{lemma}

\smallskip

\noindent{\bf Proof of Theorem \ref{teo 2}.} We argue indirectly and suppose that
$u$ is a positive supersolution of (\ref{eq 1.1-ab}). The main idea below is to reflect
$\Omega$ to a bounded punctured domain through the Kelvin transform and then to obtain a contradiction
by Lemma \ref{sub}.

Without loss of generality, we may assume that  $\Omega$ is a connected exterior domain satisfying $0\not\in \overline\Omega$.
Denote by
 $$
 \Omega^\sharp=\left\{x\in\mathbb{R}^N: \ \frac{x}{|x|^2}\in\Omega\right\}\quad{\rm and}\quad v(x)=u\left(\frac{x}{|x|^2}\right)\quad{\rm for}\  x\in\Omega^\sharp.
 $$
Clearly, $\Omega^\sharp$ is a bounded punctured domain. %Let us define  $\mathcal{D}=\Omega^\sharp\cup\{0\}$.
By direct computation, for  $x\in \Omega^\sharp$, we have
 $$
 \nabla v(x)=\nabla u\left(\frac{x}{|x|^2}\right)\frac1{|x|^2}-2\left(\nabla u\left(\frac{x}{|x|^2}\right)\cdot x\right)\frac{x}{|x|^4}
 $$
and
 $$
 \Delta v(x)=\frac1{|x|^4} \Delta u\left(\frac{x}{|x|^2}\right)
 +\frac{2(2-N)}{|x|^4} \left(\nabla u\left(\frac{x}{|x|^2}\right)\cdot x\right).
 $$

Let
 $$
 u^\sharp(x)=|x|^{2-N}v(x),\ \ \ V^\sharp(x)= |x|^{-4}V\left(\frac{x}{|x|^2}\right).
 $$
Then for $x\in \Omega^\sharp$, using the fact that $\Delta(|x|^{2-N})=0$, we observe that
 \begin{eqnarray*}
 -\Delta u^\sharp(x) &=& -\Delta v(x) \,|x|^{2-N} -2\nabla v(x)\cdot(\nabla |x|^{2-N})  \\
 &=& -\left[\frac1{|x|^4} \Delta u\left(\frac{x}{|x|^2}\right)+\frac{2(2-N)}{|x|^4}\left (\nabla u\left(\frac{x}{|x|^2}\right)\cdot x\right)\right]|x|^{2-N}\\
 &&-\frac{2(2-N)x}{|x|^2} \left[\nabla u\left(\frac{x}{|x|^2}\right)\frac1{|x|^2}-2\left(\nabla u\left(\frac{x}{|x|^2}\right)\cdot x\right)\frac{x}{|x|^4}\right]|x|^{2-N}\\
 &=& |x|^{-2-N} -\Delta u\left(\frac{x}{|x|^2}\right)\\
 &\ge& V^\sharp(x)u^\sharp(x).
 \end{eqnarray*}

Set $V^*(x)=\frac{-\Delta u^\sharp (x)}{u^\sharp(x)}$. Thus, we notice that
$V^*$ is continuous in $\Omega^\sharp$ and $V^*\ge V^\sharp$ in $\Omega^\sharp$. Therefore,
 $$
 V^*(x)\ge |x|^{-4}V(\frac{x}{|x|^2}),\quad \forall\, x\in\Omega^\sharp.
 $$
Because of (\ref{v2}), it follows that
 $$
 \liminf_{|x|\to0^+} V^*(x)|x|^2>\frac{(N-2)^2}{4},
 $$
which in turn implies that there exist $\mu_1>(N-2)^2/4$ and $r_1>0$ such that
\begin{eqnarray*}
 V^*(x)\ge \mu_1|x|^{-2},\quad \forall\, x\in B_{r_1}(0)\setminus\{0\}.
\end{eqnarray*}

Denote $u_0=u^\sharp-\varphi_0$, where $\varphi_0$ is the unique positive solution of
 \begin{eqnarray*}
 \arraycolsep=1pt\left\{
 \begin{array}{lll}
 \displaystyle -\Delta \varphi_0=0\ \
 &{\rm in}\  B_{r_1}(0) ,\\[1.5mm]
 \phantom{   -\Delta }
 \displaystyle  \varphi_0=u^\sharp \ \  &{\rm   on}\   \partial  B_{r_1}(0).
 \end{array}\right.
\end{eqnarray*}
Then $u_0$ satisfies $u_0=0$ on $\partial  B_{r_1}(0)$ and
\begin{eqnarray*}
-\Delta u_0(x)-\frac{\mu_1}{|x|^2}u_0=V^*(x)u^\sharp(x) -\frac{\mu_1}{|x|^2} u^\sharp(x)+\frac{\mu_1}{|x|^2} \varphi_0(x):=f^*(x)\geq0,\quad \forall\,  x\in B_{r_1}(0)\setminus\{0\}.
\end{eqnarray*}
As a consequence, $u_0$ is a positive solution of (\ref{eq 2.2}) with $\mathcal{D}=B_{r_1}(0)$,
$\mu=\mu_1$ and $f=f^*$. This contradicts
Lemma \ref{sub}. Thus, (\ref{eq 1.1-ab}) admits no positive supersolution, and the proof is complete. \hfill$\Box$

\subsection{Proof of Theorem \ref{coro 1}}

\noindent{\bf Proof of Theorem \ref{coro 1}.}  On the contrary, suppose that problem
\begin{eqnarray*}
 -\Delta u= \mu |x|^{-\alpha} u\quad{\rm  in}\ \mathbb{R}^N\setminus\{0\}
\end{eqnarray*}
has a positive supersolution $u_0$, that is,
\begin{equation}\label{2.27}
 -\Delta u_0\ge  \mu |x|^{-\alpha} u_0\quad{\rm  in}\ \mathbb{R}^N\setminus\{0\}\ {\rm pointwise}.
\end{equation}

When $\alpha=2$ and $\mu>\frac{(N-2)^2}{4}$, a contradiction can be seen directly from  Theorem \ref{teo 2}.

When $\alpha< 2$, (\ref{2.27}) can be written as
\begin{eqnarray*}
 -\Delta u_0 \ge \mu |x|^{2-\alpha} |x|^{-2} u_0\ge(N-2)^2 |x|^{-2} u_0 \quad{\rm  in}\   \mathbb{R}^N\setminus B_{r_\mu}(0),
\end{eqnarray*}
where
 $$
 r_\mu=\left[\frac{(N-2)^2}{\mu}\right]^{\frac1{2-\alpha}}>0.
 $$
This contradicts Theorem \ref{teo 2}.

When $\alpha>2$, it follows from (\ref{2.27})  that there exists $\tilde f\ge 0$ such that
\begin{eqnarray*}
 -\Delta u_0 &=&   \mu |x|^{2-\alpha} |x|^{-2} u_0 +\tilde f
  \\  &\ge&   (N-2)^2 |x|^{-2} u_0\ \quad{\rm  in}\ B_{r_\mu}(0)\setminus\{0\}.
\end{eqnarray*}
Denote by $\varphi_0$ the unique positive solution of
 \begin{eqnarray*}
 \arraycolsep=1pt\left\{
 \begin{array}{lll}
 \displaystyle -\Delta \varphi_0=0\ \
 &{\rm in}\  B_{r_\mu}(0) ,\\[1.5mm]
 \phantom{   -\Delta }
 \displaystyle  \varphi_0=u_0 \ \  &{\rm   on}\   \partial  B_{r_\mu}(0).
 \end{array}\right.
\end{eqnarray*}
Then $\varphi_0\in C^{2}(B_{r_\mu}(0))\cap C(\overline B_{r_\mu}(0))$, and
$v_0:=u_0-\varphi_0$ is  bounded from below and  is a solution of
 \begin{eqnarray*}
 \arraycolsep=1pt\left\{
 \begin{array}{lll}
 \displaystyle -\Delta v_0=\overline f\ \ \ \
 &{\rm in}\  B_{r_\mu}(0)\setminus \{0\},\\[1.5mm]
 \phantom{   -\Delta }
 \displaystyle  v_0=0 \ \  &{\rm   on}\   \partial B_{r_\mu}(0),
 \end{array}\right.
 \end{eqnarray*}
where $\overline f= \mu|x|^{-\alpha} u_0+\tilde f>0$.

Next we show $v_0>0$ in $B_{r_\mu}(0)\setminus \{0\}$. Indeed, consider
 \begin{eqnarray*}
 \arraycolsep=1pt\left\{
 \begin{array}{lll}
 \displaystyle -\Delta z=\mu|x|^{-\alpha} u_0 \chi_{_{B_{r_\mu}(0)\setminus B_{{r_\mu}/2}(0)}}\ \ \ \
 &{\rm in}\  B_{r_\mu}(0),\\[1.5mm]
 \phantom{   -\Delta }
 \displaystyle  z=0 \ \  &{\rm   on}\   \partial B_{r_\mu}(0),
 \end{array}\right.
 \end{eqnarray*}
which admits a unique positive bounded solution, denoted by $z_0$.
For any small $\epsilon>0$, set
 $$
 \psi_\epsilon(x)=\epsilon(r_\mu^{2-N}-|x|^{2-N})+z_0(x).
 $$
Then there exists $r_\epsilon\in(0,r_\mu)$ such that
$\lim_{\epsilon\to0^+}r_\epsilon=0$ and
$\psi_\epsilon\le v_0$ on $\overline B_{r_\epsilon}(0)$. By the classical comparison principle, we have
 $$
 v_0\ge \psi_\epsilon  \quad{\rm in}\   B_{r_\mu}(0)\setminus{B_{r_\epsilon}(0)}.
 $$
Sending $\epsilon\to0$ , we can conclude $v_0\ge z_0>0$ in $B_{r_\mu}(0)\setminus \{0\}$.

Furthermore,  we see that  $v_0=0$  on $\partial B_{r_\mu}(0)$ and
 \begin{eqnarray*}
 -\Delta v_0  &= & \mu |x|^{-\alpha}  u_0+\tilde f \\
 &\geq&  (N-2)^2 |x|^{-2}v_0+F\ \quad{\rm  in}\  B_{r_\mu}(0)\setminus\{0\},
 \end{eqnarray*}
where
 $$
 F(x)=\mu |x|^{-\alpha}\varphi_0+ \tilde f +(\mu |x|^{-\alpha} -(N-2)^2 |x|^{-2}) u_0\ge0 \quad {\rm for}\ 0<|x|<r_\mu.
 $$
 Thus, we obtain a contradiction with Theorem \ref{sub}  and the proof ends. \hfill$\Box$

\begin{remark} With a slight modification of the proof for the case  $\alpha>2$ in Theorem \ref{coro 1},
we can observe the following result: for $\mu>(N-2)^2/4$ with the dimension $N\ge 2$, and $f\in L^\infty_{loc}(\overline{\mathcal{D}\setminus \{0\}})$ a nonnegative function with
$\mathcal{D}$ a bounded  smooth domain containing the origin, then the Hardy problem (\ref{eq 2.2})
has no nontrivial nonnegative supersolution.

\end{remark}

\setcounter{equation}{0}
\section{Nonexistence of positive supersolution of nonlinear problems }

\subsection{Proof of Theorem \ref{teo 1}}

Our proof is based on Theorem \ref{teo 2} and the following comparison principle.

\begin{lemma}\label{lm cp}%\cite[Lemma 2.1]{CQZ}
Assume that  $\Omega$  is a smooth exterior domain, $f_1,\,f_2$ are continuous in $\Omega$,
$g_1,\,g_2$ are continuous on $\partial \Omega$, and
$f_1\ge f_2\ {\rm in}\ \Omega$ and $g_1\ge g_2\ {\rm on}\ \partial \Omega.$
Let $u_1$ and $u_2$ satisfy $-\Delta  u_1\geq f_1$ in $\Omega$, $u_1\geq g_1$ on $\partial{\Omega}$, and
$-\Delta  u_2\leq f_2$ in $\Omega$, $u_2\leq g_2$ on $\partial{\Omega}$. If $\liminf_{|x|\to+\infty}u_1(x) \ge \limsup_{|x|\to+\infty}u_2(x),$ then we have
$u_1\ge u_2\  {\rm on}\ \overline\Omega$.

\end{lemma}
 {\bf Proof.}  Such a comparison principle may be folklore; we here provide a simple proof for sake of completeness.
Letting $w=u_2-u_1$, then $w$ satisfies
 $$
 -\Delta w \le 0\ \  {\rm in}\ {\Omega}, \ \ w\le 0 \ \ {\rm   on}\  \partial{\Omega}\ \ \mbox{and}\ \
 \limsup_{|x|\to+\infty}w(x) \le0.
 $$
Thus, given $\epsilon>0$, there exists $r_\epsilon>0$ converging to infinity  as $\epsilon\to0$ such that
 $$
 w\le \epsilon (r_\epsilon^{2-N} +1)\quad{\rm on}\ \partial B_{r_\epsilon}(0).
 $$
As a result, we have
 $$
 w(x)\le \epsilon< \epsilon (|x|^{2-N} +1) \quad{\rm on}\  \partial(\Omega\cap B_{r_\epsilon}(0)).
 $$
Note that $\Delta(|x|^{2-N})=0$. It then follows from the classical comparison principle in any bounded smooth domain that
 $$
 w(x)\le \epsilon (|x|^{2-N} +1)\quad{\rm in}\ \Omega \cap B_{r_\epsilon}(0).
 $$
According to the arbitrariness of $\epsilon>0$, we can conclude that $w\le 0$ on $\overline\Omega$.
\hfill$\Box$\medskip

\vskip10pt
In order to prove Theorem \ref{teo 1}, we introduce the auxiliary function
 $$
 \tilde w_0(r)= r^{2-N} (\ln  r)^{\sigma}\ \ \ \mbox{with\ $\sigma>0$}.
 $$
Elementary computation yields
 $$
 \tilde w_0'(r)=(2-N)r^{1-N} (\ln r)^\sigma+\sigma r^{1-N} (\ln r)^{\sigma-1}
 $$
and
 $$
 \tilde w_0''(r)=(2-N)(1-N)r^{-N} (\ln r)^\sigma+\sigma(3-2N) r^{-N} (\ln r)^{\sigma-1}+\sigma(\sigma-1) r^{-N} (\ln r)^{\sigma-2}.
 $$
Let $w_0(x)=\tilde w_0(|x|)$. Then, if $x$ satisfies
$|x|>\max\Big\{1,\  e^{\frac{2(1-\sigma)}{(N-2)}}\Big\},$
we have
 $$
 \nabla w_0(x)=\tilde w_0'(|x|)\nabla|x|=\tilde w_0'(|x|)\frac{x}{|x|},
 $$
and so
 \begin{eqnarray}
 -\Delta w_0(x)
 &=&-\Big[w_0''(|x|)+\tilde w_0'(|x|)\cdot\frac{N-1}{|x|}\Big] \nonumber \\
 &=& \sigma(N-2)|x|^{-N} (\ln |x|)^{\sigma-1}-\sigma(\sigma-1) |x|^{-N} (\ln |x|)^{\sigma-2} \nonumber \\
 &\le &  \frac32\sigma(N-2)|x|^{-N} (\ln |x|)^{\sigma-1}.
 \label{3.0}
 \end{eqnarray}

Now, with the aid of the function $w_0$, we are ready to prove Theorem \ref{teo 1}.

\vskip10pt

\noindent{\bf Proof of Theorem \ref{teo 1}.} We use a contradiction argument and
suppose that (\ref{eq 1.1}) has a positive supersolution $u$. Since
$h(x,u)\ge \tilde h (x,u)$ for  $(x,u)\in (\mathbb{R}^N\setminus B_e(0))\times [0,+\infty)$, then $u$ fulfills
 \begin{equation}\label{eq 4.1}
 -\Delta u\ge \tilde h (x,u)\quad{\rm in}\ \, \mathbb{R}^N\setminus B_e(0).
 \end{equation}
By the positivity of $u$, one can find a constant $c_0>0$ such that
 $$
 u(x)\ge c_0e^{2-N},\quad  \forall\, x\in \partial B_e(0).
 $$
Hence, Lemma \ref{lm cp} gives
\begin{equation}\label{eq 4.1a}
 u(x)\ge c_0|x|^{2-N},\quad \forall\, |x|\ge e.
 \end{equation}

By our assumption (H)-(a), one can see that
 \begin{eqnarray*}
 \frac{ \tilde h  (x,u(x))}{u(x)} & \ge &
 \frac{\tilde h  (x,c_0|x|^{2-N}) }{c_0|x|^{2-N}},\quad \forall\, |x|\ge e,
 \end{eqnarray*}
and
if (H)-(b) holds, then
 \begin{equation}\label{4.2}
 \liminf_{|x|\to+\infty}\frac{\tilde h  (x,c_0|x|^{2-N}) }{c_0|x|^{2-N}}|x|^2>\frac{(N-2)^2}{4}.
 \end{equation}
Take  $V(x)=\frac{\tilde h (x,u(x))}{u(x)}$, which satisfies $\liminf_{|x|\to+\infty} V(x)|x|^2>(N-2)^2/4$.
It is easily observed that $u$ is a positive supersolution of
 $$
 -\Delta u=V(x) u \quad   {\rm in}\ \mathbb{R}^N\setminus  B_e(0).
 $$
This contradicts Theorem \ref{teo 2}. Thus, (\ref{eq 4.1}) has no positive
supersolution provided (H) holds.

If (H)-(b) fails, in the sequel we shall establish the nonexistence result in
Theorem \ref{teo 1} using the assumption (H)-(b1)(b2).
Under (H)-(b1), by (\ref{eq 4.1a}), there exists $\varrho_0\ge e$ such that
 \begin{eqnarray}\label{4.2a}
 \tilde h  (x,u(x)) & \ge & \tilde h  (x,c_0|x|^{2-N})\frac{u(x)}{c_0|x|^{2-N}}\nonumber\\
 &\ge&\tilde h  (x,c_0|x|^{2-N})\nonumber \\
 &\ge & m_0|x|^{-N} (\ln  |x|)^{-\sigma_1} ,\quad\ \forall\, |x|\ge \varrho_0,
 \end{eqnarray}
where
 $$
 m_0=\min\Big\{1,\  \liminf_{|x|\to+\infty} \tilde h (x,\,c_0|x|^{2-N})  |x|^{N} (\ln|x|)^{\sigma_1}\Big\}>0.
 $$
In light of (\ref{3.0}) (by taking $\sigma=1-\sigma_1>0$ there),
one can find $t_1\in(0, \frac{m_0}{3(N-2)\sigma_1})$ and $\varrho_1\geq\varrho_0$ such that
 $$
 u(x)\ge t_1w_0(x),\quad |x|=\varrho_1
 $$
and
 $$
 -\Delta u(x)\ge m_0|x|^{-N} (\ln  |x|)^{-\sigma_1}\ge -\Delta (t_1w_0(x)) ,\quad \forall\, |x|>\varrho_1.
 $$
Then by Lemma \ref{lm cp}, it follows that
 \begin{equation}\label{4.1}
 u(x)\ge t_1|x|^{2-N}(\ln|x|)^{\theta_1}, \quad \forall\, |x|\ge \varrho_1
 \end{equation}
with $\theta_1=1-\sigma_1$.

As a next step, we are going to improve the decay of $u$ at infinity by an induction argument. To this end, let $\{\theta_j\}_j$ be the sequence generated by
 \begin{equation}\label{5.1}
 \theta_{j+1}=\theta_j+\theta_1=(j+1)\theta_1, \quad  j=1,2,3\cdots.
 \end{equation}
In view of (H)-(a), we may assume that $\sigma_2\geq1$, where
$\sigma_2>0$ appears from the assumption (H)-(b2). Furthermore, by observing that $\lim_{j\to+\infty} \theta_j=+\infty$,
one can assert that there exists $j_0\in\mathbb{N} $  such that
 \begin{equation}\label{2.3}
 \theta_{j_0} \ge \sigma_2\quad {\rm and}\quad \theta_{j_0-1}<\sigma_2.
 \end{equation}

We now claim that for any integer $j$, there exist $\varrho_j\geq\varrho_0$ and $t_j>0$ such that
 \begin{equation}\label{5.2j}
 u(x)\ge t_j|x|^{2-N}(\ln|x|)^{\theta_j}, \quad \forall\,|x|\ge \varrho_j.
 \end{equation}
Indeed, when $j=1$, (\ref{5.2j}) has been proved above.
Assume that (\ref{5.2j}) holds for some $j$, we will show that (\ref{5.2j}) holds true for $j+1$.
Notice that there exists $r_{j}\ge \varrho_j$ such that
 $$
 t_j|x|^{2-N}(\ln|x|)^{\theta_j} \ge c_0 |x|^{2-N}, \quad \forall\, |x|\ge r_{j}.
 $$
By (H)-(a)(b1), (\ref{4.2a}) and (\ref{5.2j}), we obtain
\begin{eqnarray*}
\tilde h  (x,u(x)) & \ge & \tilde h  (x,\, t_j|x|^{2-N}(\ln|x|)^{\theta_j})\frac{u(x)}{t_j|x|^{2-N}(\ln|x|)^{\theta_j}}
\\&\ge&\tilde h  (x,\, t_j|x|^{2-N}(\ln|x|)^{\theta_j})
\\&\ge& \tilde h  (x,\, c_0|x|^{2-N})\frac{t_j|x|^{2-N}(\ln|x|)^{\theta_j}}{c_0|x|^{2-N}}
\\&\ge& m_0\frac{t_j}{c_0} |x|^{-N}(\ln|x|)^{\theta_j-\sigma_1},\ \ \forall\, |x|\ge r_{j}.
\end{eqnarray*}

Taking $\sigma=1+\theta_j-\sigma_1>0$ in (\ref{3.0}), and then using
Lemma \ref{lm cp} (by comparing $u$ with $|x|^{2-N}(\ln|x|)^{\theta_j+\theta_1}$),
we can conclude that
\begin{eqnarray*}
 u(x)\ge t_{j+1}|x|^{2-N}(\ln|x|)^{\theta_j+\theta_1}, \quad \forall\,|x|\ge \varrho_{j+1},
\end{eqnarray*}
for some $t_{j+1}>0$ and $\varrho_{j+1}\ge  r_{j} $.
This verifies the previous claim (\ref{5.2j}).

Therefore,  (\ref{2.3}) and (\ref{5.2j}) imply that
 $$
 u(x)\ge t|x|^{2-N}(\ln|x|)^{\sigma_2}, \quad \forall\, |x|\ge \varrho_{j_0},
 $$
and in turn by (H)-(a),
 \begin{equation}\label{3.2}
 \frac{\tilde h  (x,u(x))}{u(x)}\ge\frac{\tilde h(x,t|x|^{2-N}(\ln|x|)^{\sigma_2} )}{t|x|^{2-N}(\ln|x|)^{\sigma_2}},\quad\forall\, x\in\mathbb{R}^N\setminus  B_{\varrho_{j_0}}(0).
\end{equation}
By (H)-(b2), there exist $\varrho>\varrho_{j_0}$ and $\epsilon_0>0$  such that
 $$
 \frac{\tilde h  (x,u(x))}{u(x)}\ge \Big(\frac{(N-2)^2}4+\epsilon_0\Big)\frac1{|x|^2}, \quad \forall\,  x\in\mathbb{R}^N\setminus  B_{\varrho}(0).
 $$
Let  $V(x)=\frac{h(x,u(x))}{u(x)}$. Then
 $$
 V(x)|x|^2>\frac{(N-2)^2}{4}+\epsilon_0, \quad \forall\,  x\in\mathbb{R}^N\setminus  B_{\varrho}(0).
 $$
Clearly, $u$ is a positive supersolution of
 $$-\Delta u=V(x) u \quad   {\rm in}\  \mathbb{R}^N\setminus  B_{\varrho}(0),
 $$
which is a contradiction with Theorem \ref{teo 2}.
Thus,  (\ref{eq 1.1}) admits no positive  supersolution. The proof is complete.\hfill$\Box$

\smallskip

\subsection{Proof of Theorem \ref{teo 3}}

\noindent{\bf Proof of Theorem \ref{teo 3}.} Suppose that (\ref{sys 1.1})
has a positive supersolution $(u,v)$. Clearly, $(u,v)$ fulfills
 \begin{equation}\label{eq 4.2}
 -\Delta u\ge \tilde h_1 (x,u,v)\geq0\quad {\rm and} \quad -\Delta v\ge \tilde h_2 (x,u,v)\geq0
 \quad{\rm in} \ \mathbb{R}^N\setminus B_e(0).
 \end{equation}
Since one can find a small $c_0>0$ such that
 $$
 u(x),\ v(x)\ge c_0e^{2-N},\quad  \forall\, x\in \partial B_e(0),
 $$
an analysis similar to that of obtaining (\ref{eq 4.1a}) yields that
 \begin{equation}\label{eq 4.2-a}
 u(x),\ v(x)\ge c_0|x|^{2-N},\quad \forall\, |x|\ge e.
 \end{equation}
Thus, it follows from (SH)-(a) that
 \begin{eqnarray*}
 \frac{ \tilde h_1  (x,\,u(x),\,v(x))}{u(x)}\ge \frac{ \tilde h_1  (x,\,u(x),\,c_0|x|^{2-N})}{u(x)}\ge\frac{\tilde h_1  (x,\,c_0|x|^{2-N},\,c_0|x|^{2-N} ) }{c_0|x|^{2-N}},\ \ \forall\, |x|\ge e
 \end{eqnarray*}
and
 \begin{eqnarray*}
 \frac{ \tilde h_2  (x,\,u(x),\,v(x))}{v(x)}\ge \frac{ \tilde h_2  (x,\,c_0|x|^{2-N},\,v(x))}{v(x)}   \ge   \frac{\tilde h_2  (x,\,c_0|x|^{2-N},\,c_0|x|^{2-N} ) }{c_0|x|^{2-N}},\ \ \forall\, |x|\ge e.
 \end{eqnarray*}

If (SH)-(b) holds, we then have
 \begin{equation}\label{4.2-1}
 \liminf_{|x|\to+\infty}\frac{\tilde h_1  (x,c_0|x|^{2-N}, c_0|x|^{2-N}) }{c_0|x|^{2-N}}|x|^2>\frac{(N-2)^2}{4}
 \end{equation}
or
 \begin{equation}\label{4.2-2}
 \liminf_{|x|\to+\infty}\frac{\tilde h_2  (x,c_0|x|^{2-N},c_0|x|^{2-N}) }{c_0|x|^{2-N}}|x|^2>\frac{(N-2)^2}{4}.
 \end{equation}
By taking  $V(x)=\frac{\tilde h_1 (x,\,u(x),\,v(x))}{u(x)}$ if (\ref{4.2-1}) holds (or $V(x)=\frac{\tilde h_1 (x,u(x),\,v(x))}{v(x)}$ if (\ref{4.2-2}) holds), we see that $u$ (or $v$) is a positive supersolution of
 $$
 -\Delta u=V(x) u \quad   {\rm in}\  \mathbb{R}^N\setminus  B_e(0)
 $$
with $\liminf_{|x|\to+\infty} V(x)|x|^2>\frac{(N-2)^2}{4}$. This is impossible due to Theorem \ref{teo 2}.
Therefore,  (\ref{sys 1.1}) has no positive  supersolution.

If (SH)-(b) fails, we continue to prove the nonexistence result in Theorem  \ref{teo 3} under the assumption (SH)-(b1)(b2)(b3).
There are three cases to distinguish as follows.

{\it Case 1: $0<\sigma_3<1,\ \sigma_4\geq1$.} By (\ref{eq 4.2-a}) and the assumption (SH)-(a), we have
 \begin{eqnarray*}
 \tilde h_1  (x,\, u(x),\, v(x))\ge\tilde h_1  (x,\, u(x),\,c_0|x|^{2-N}):=\tilde  h(x,u),\ \ \forall\, |x|>e.
 \end{eqnarray*}
Then  $u$ verifies that
 \begin{equation}\label{eq 3.1}
 -\Delta u\ge \tilde h(x,u) \quad   {\rm in}\  \,\mathbb{R}^N\setminus  B_e(0).
 \end{equation}
By virtue of the assumption (SH)-(a)(b1), there exists $\rho_0\ge e$ such that
 \begin{eqnarray}
 \tilde h(x,u) % & \ge & \tilde h_1  (x,\,c_0|x|^{2-N},\,c_0|x|^{2-N})\frac{u(x)}{c_0|x|^{2-N}} \nonumber \\
 &\ge&\tilde h_1  (x,\,c_0|x|^{2-N},\,c_0|x|^{2-N}) \label{3.1}\\
 &\ge & \frac12 m_1|x|^{-N}    (\ln |x|)^{-\sigma_3}  ,\quad\ \forall\,|x|\ge \rho_0 ,\nonumber
 \end{eqnarray}
where
 $$
 m_1=\min\Big\{1,\ \liminf_{|x|\to+\infty} \tilde h_1  (x,\,c_0|x|^{2-N},\,c_0|x|^{2-N})  |x|^{N}  (\ln |x|)^{\sigma_{3}}\Big\}.
 $$
This implies that (H)-(b1) holds. Moreover,  (SH)-(b2) indicates that $\tilde h$ satisfies (H)-(b2).
Thus, an application of Theorem \ref{teo 1} to problem (\ref{eq 3.1}) leads to a contradiction. Hence,
(\ref{sys 1.1}) has no positive  supersolution in Case 1.

{\it Case 2: $\sigma_3\geq1,\ 0<\sigma_4<1$.} The proof is similar to  {\it Case 1. }

{\it Case 3: $0<\sigma_3,\,\sigma_4<1$.} Due to (SH)-(b1), we can deduce
 \begin{eqnarray*}
 \tilde h_i (x,\, u(x),\, v(x))  \ge \frac12 m_i|x|^{-N}   (\ln  |x|)^{-\sigma_{2+i}},\quad \forall\,|x|\ge \rho_0,
 \end{eqnarray*}
for some $\rho_0>e$ and
 $$
 m_i=\min\Big\{1,\ \liminf_{|x|\to+\infty} \tilde h_i  (x,\,c_0|x|^{2-N}, c_0|x|^{2-N})  |x|^{N}  (\ln |x|)^{\sigma_{2+i}}\Big\}.
 $$
Proceeding similarly as in (\ref{4.1}), one can assert that for some $t_0>0$ and $\rho_1>\rho_0$,
 \begin{equation}\label{4.1-s}
 u(x)\ge t_0|x|^{2-N} (\ln|x|)^{1-\sigma_3} \quad{\rm and} \quad v(x)\ge t_0|x|^{2-N} (\ln|x|)^{1-\sigma_4}
 \quad {\rm in}\ \mathbb{R}^N\setminus  B_{\rho_1}(0).
 \end{equation}
Then, reasoning as in the part of the claim (\ref{5.2j}), we have
 $$
 u(x)\ge t_j|x|^{2-N}(\ln |x|)^{ j(1-\sigma_3)}\quad{\rm and}\quad v(x)
 \ge t_j|x|^{2-N}(\ln |x|)^{j(1-\sigma_4)},\quad \forall\,|x|\ge \rho_j
 $$
for a sequence $\{(t_j,\rho_j)\}_{j=1}^\infty$. Therefore, there exists a large integer $j^*$ such that
 $$
 j^*(1-\sigma_3)>\sigma_5\quad{\rm and}\quad j^*(1-\sigma_4)>\sigma_6.
 $$
Thus, we obtain
 $$
 u(x)\ge t_{j^*}|x|^{2-N}(\ln |x|)^{ \sigma_5},\quad \ v(x)\ge t_{j^*}|x|^{2-N}(\ln |x|)^{ \sigma_6},
 \quad \forall\,|x|\ge \rho_{j^*}.
 $$
This then yields
 $$
 \frac{\tilde h_1  (x,u(x),v(x))}{u(x)}\ge\frac{\tilde h_1(x,t_{j^*}|x|^{2-N}(\ln|x|)^{\sigma_5},t_{j^*}|x|^{2-N}(\ln|x|)^{\sigma_6}) }{t_{j^*}|x|^{2-N}(\ln|x|)^{\sigma_5}},
 \quad \forall\,|x|\ge \rho_{j^*}
 $$
and
 $$
 \frac{\tilde h_2  (x,u(x),v(x))}{u(x)}\ge\frac{\tilde h_2(x,t_{j^*}|x|^{2-N}(\ln|x|)^{\sigma_5},t_{j^*}|x|^{2-N}(\ln|x|)^{\sigma_6}) }{t_{j^*}|x|^{2-N}(\ln|x|)^{\sigma_5}},
 \quad \forall\,|x|\ge \rho_{j^*}.
 $$

By the assumption (SH)-(b3), there exist $\varrho>\rho_{j^*}$ and $\epsilon_0>0$ such that
 $$
 \frac{\tilde h_1  (x,u(x),v(x))}{u(x)}\,\left(\ {\rm or} \ \ \frac{\tilde h_2  (x,u(x),v(x))}{v(x)}\right)>(\frac{(N-2)^2}4+\epsilon_0)\frac1{|x|^2},
 \quad   \forall\, x\in\mathbb{R}^N\setminus  B_{\varrho}(0).
 $$
Taking  $V(x)=\frac{h_1(x,u(x),v(x))}{u(x)}$ (or $V(x)=\frac{h_2(x,u(x),v(x))}{v(x)}$), $u$ is a positive supersolution of
 $$
 -\Delta u=V(x) u \quad   {\rm in}\ \mathbb{R}^N\setminus  B_{\varrho}(0),
 $$
contradicting Theorem \ref{teo 2}. As a consequence,  (\ref{sys 1.1}) has no positive  supersolution in Case 3.
The proof is now complete. \hfill$\Box$

\setcounter{equation}{0}
\section{Application: two examples}

In this section, we shall use two typical examples to illustrate the optimality of the nonexistence results obtained
by this paper.\medskip

{\bf Example 1:} $h(x,u)=|x|^{\beta} (\ln |x|)^\tau u^{p}$.

\medskip When $\beta>-2$,  then $p_\beta^*=\frac{N+\beta}{N-2}>1$. We have

\begin{proposition}\label{cr 2} Assume that  $h(x,u)=|x|^{\beta} (\ln |x|)^\tau u^{p},\,x\in \mathbb{R}^N\setminus B_{e}(0)$ with $\beta>-2$.
Then the following assertions hold.

 \begin{itemize}
\item[$\rm (i)$] Problem (\ref{eq 1.1})  has no positive supersolution provided that either
 $1\le p<p_\beta^*, \,\tau\in\mathbb{R}$ or $p=p_\beta^*, \,\tau>-1;$

 \item[$\rm (ii)$] Problem (\ref{eq 1.1})  has no positive bounded supersolution provided that
 $p\in(-\infty,1),\,\tau\in\mathbb{R}$.
 \end{itemize}

\end{proposition}

\noindent{\bf Proof.} We shall apply Theorem \ref{teo 1} to obtain the desired results.
It suffices to check the condition (H).

We first verify (i). When $p=1$, the nonexistence follows from  Theorem \ref{teo 2} directly.
Clearly, the condition (H)-(a) is fulfilled once $p>1$. We also note that
 $$
 h(x,|x|^{2-N})|x|^N=|x|^{\beta+N-(N-2)p}(\ln |x|)^\tau,\quad \forall\, |x|>e.
 $$
If $p<p^*_\beta$, then $\beta+N-(N-2)p>0$. So for any $t>0$ and $\tau\in\mathbb{R}$, then
 $$
 \lim_{|x|\to+\infty}  h  (x,\,t|x|^{2-N})  |x|^{N} =+\infty.
 $$
Thus, the assumption (H)-(b) is satisfied.

If $p=p^*_\beta$ and $\tau\in(-1,0)$, by taking $\sigma_1=-\tau>0$, we have
 $$
 h (x,\,t|x|^{2-N})  |x|^{N}  (\ln|x|)^{\sigma_1}=t^{p_\beta^*}(\ln|x|)^{\tau+\sigma_1}=t^{p_\beta^*},\quad \forall\, |x|>e.
 $$
If $p=p^*_\beta$ and $\tau\ge0$, by taking $\sigma_1=\frac12>0$, we have
 $$
 h (x,\,t|x|^{2-N})  |x|^{N}  (\ln|x|)^{\sigma_1}=t^{p_\beta^*}(\ln|x|)^{\tau+\frac12}\ge t^{p_\beta^*},\quad \forall\, |x|>e.
 $$

Furthermore, let us choose $\sigma_2>0$ such that $\tau +\sigma_2(p_\beta^*-1)>0$. Then
 $$
 \frac{h(x,\,t|x|^{2-N}(\ln|x|)^{\sigma_2})}{t|x|^{-N}(\ln|x|)^{\sigma_2}}
 =t^{p_\beta^*-1}(\ln|x|)^{\tau +\sigma_2(p_\beta^*-1)}, \quad \forall\, |x|>e.
 $$
Therefore, $h$ satisfies the assumption (H)-(b1)(b2). Thus, the assertion (i) is proved.

We next verify the assertion (ii). Arguing indirectly, we suppose that (\ref{eq 1.1})
has a positive bounded supersolution $u$. Denote
 $$
 \tilde h(x,u)= M^{p-1} |x|^{\beta}(\ln|x|)^\tau u,\ \ \mbox{if}\ p\in(-\infty,1],\,\tau\in\mathbb{R},
 $$
where $M=1+\sup_{x\in\Omega} u(x)>1$.

Then, it is easily checked that
 $$
 h(x,u)=|x|^{\beta} (\ln |x|)^\tau u^{p}\geq M^{p-1} |x|^{\beta}(\ln|x|)^\tau u
 =\tilde h(x,u),\quad \forall\, |x|>e,\ u\in(0,M],
 $$
and $\tilde h$ fulfills the assumption (H)-(a)(b) for any $\tau\in\mathbb{R}$.
Consequently, $u$ is a positive supersolution of
 \begin{equation}\label{eq 2.3}
 -\Delta u=\tilde h(x,u)\quad {\rm in}\ \,\mathbb{R}^N\setminus B_{e}(0).
 \end{equation}
However, (\ref{eq 2.3}) has no positive supersolution by Theorem \ref{teo 1}.
Such a contradiction implies that (\ref{eq 1.1}) has no positive bounded supersolution.
The assertion (ii) follows. \hfill$\Box$

\bigskip

\begin{proposition}\label{cr 2-5}
 Assume that  $h(x,u)=|x|^{\beta} (\ln |x|)^\tau u^{p}, \,\forall\, x\in \mathbb{R}^N\setminus B_\ell(0)$, where $\beta\in\mathbb{R}$.
Then for some large $\ell>e$, problem (\ref{eq 1.1}) with $\Omega=\mathbb{R}^N\setminus B_\ell(0)$ has a positive bounded supersolution if one of the following conditions is satisfied:
 \begin{itemize}
 \item[$\rm (i)$]   $p> p^*_\beta$ and $\tau\in\mathbb{R}$;
 \item[$\rm (ii)$]   $p= p^*_\beta$ and $\tau<-1$.

%  \item[$\rm (ii)$] $\beta\leq-2$, $p\geq p^*_\beta$ and either $p\geq1,\,\tau<-1$ or $p<1,\,\tau\in\mathbb{R}$.
 \end{itemize}
\end{proposition}

\noindent{\bf Proof.} Recall that $w_0(x)=|x|^{2-N} (\ln|x|)^\sigma$ with $\sigma>0$.
If
 $$
 |x|>\max\Big\{1,\  e^{\frac{2(\sigma-1)}{N-2}}\Big\},
 $$
it immediately follows from (\ref{3.0}) that
\begin{eqnarray*}
 -\Delta w_0(x) &=&  \sigma(N-2)|x|^{-N} (\ln |x|)^{\sigma-1}-t\sigma(\sigma-1) |x|^{-N} (\ln |x|)^{\sigma-2} \nonumber \\
   &\ge &  \frac12 \sigma (N-2)|x|^{-N} (\ln |x|)^{\sigma-1}
\end{eqnarray*}
and
\begin{eqnarray*}
h(x,w_0) &=& |x|^{\beta} (\ln |x|)^\tau  (|x|^{2-N} (\ln|x|)^\sigma)^{p}= |x|^{\beta+(2-N)p} (\ln |x|)^{\tau+\sigma p}.%  \\
  % &\leq&   t^{p}|x|^{-N} (\ln |x|)^{\tau+\sigma p},
\end{eqnarray*}

When  $p> p^*_\beta$, then $\beta+(2-N)p<-N$. One can easily see that there exists a large constant $\ell>e$ such that
 $$
 -\Delta w_0(x)\ge h(x,w_0(x)),\quad \forall\, |x|>\ell.
 $$
Hence, $w_0$ is a desired supersolution.

When  $p= p^*_\beta$, we have $\beta+(2-N)p=-N$. Similarly as above, for some large $\ell>e$,
problem (\ref{eq 1.1}) has a positive supersolution $w_0$ if $\sigma-1>\tau+\sigma p^*_\beta$, that is,
\begin{equation}\label{6.1}
 \sigma(p^*_\beta-1)<-\tau-1,
\end{equation}
where $-\tau-1>0$ by our assumption $\tau<-1$. So when $\beta>-2$, (\ref{6.1}) holds if
we take $\sigma>\frac{-\tau-1}{p^*_\beta-1}$. When $\beta\le -2$, then $p^*_\beta-1\le 0$ and
(\ref{6.1}) is satisfied if we just take $\sigma=1$.

In each case, the supersolution $w_0$ is bounded in $\mathbb{R}^N\setminus B_{\ell}(0)$.
 \hfill$\Box$ \bigskip

Our second example is the following one:

\smallskip

{\bf Example 2:} $h_1(x,u,v)=|x|^{\beta_1} (\ln |x|)^{\tau_1}u^{p_1}v^{q_1}$, \  $h_2(x,u,v)=|x|^{\beta_2} (\ln |x|)^{\tau_2}u^{p_2}v^{q_2}$.

\begin{proposition}\label{cr 2-3} Assume that
 $$
 h_1(x,u,v)=|x|^{\beta_1} (\ln |x|)^{\tau_1}u^{p_1}v^{q_1}, \quad \
 h_2(x,u,v)=|x|^{\beta_2} (\ln |x|)^{\tau_2}u^{p_2}v^{q_2},\quad \forall\,\, x\in \mathbb{R}^N\setminus B_{e}(0),
 $$
where $p_1,\, q_2\ge1$, $p_2,\, q_1\ge0$, $\beta_1,\, \beta_2>-2$ and $\tau_1,\, \tau_2\in\mathbb{R}$.
Problem (\ref{sys 1.1})  has no positive supersolution if one of the following conditions holds:

 \begin{itemize}
 \item[$\rm (i)$] $p_1+q_1<p_{\beta_1}^*,\ \tau_1,\tau_2\in\mathbb{R};$

  \item[$\rm (ii)$] $p_2+q_2<p_{\beta_2}^*,\ \tau_1,\tau_2\in\mathbb{R};$

  \item[$\rm (iii)$] $p_1+q_1=p_{\beta_1}^*,\, p_2+q_2=p_{\beta_2}^*$, $p_1>1,\,\tau_1>-1,\tau_2\in\mathbb{R};$

  \item[$\rm (iv)$] $p_1+q_1=p_{\beta_1}^*,\, p_2+q_2=p_{\beta_2}^*$, $q_2>1,\,\tau_2>-1,\tau_1\in\mathbb{R};$

  \item[$\rm (v)$] $p_1+q_1=p_{\beta_1}^*,\, p_2+q_2=p_{\beta_2}^*,%$, $p_1=q_2=1, p_2+q_1>0,
  \tau_1>-1,\, \tau_2>-1$.
 \end{itemize}

\end{proposition}
{\bf Proof.} In order to apply Theorem \ref{teo 3}, we only need to check that the nonlinearities
 $h_1,\, h_2$ satisfy (SH).

First of all, when $p_1,\, q_2\ge1$, $h_1,\, h_2$ satisfy (SH)-(a). We further note that
 $$
 h_1(x,s|x|^{2-N},s|x|^{2-N})|x|^N=s^{p_1+q_1}|x|^{\beta_1+N-(N-2)(p_1+q_1)}(\ln |x|)^{\tau_1},\quad\forall\, |x|>e,
 $$
 $$
 h_2(x,s|x|^{2-N},s|x|^{2-N})|x|^N=s^{p_2+q_2}|x|^{\beta_2+N-(N-2)(p_2+q_2)}(\ln |x|)^{\tau_2},\quad\forall\, |x|>e.
 $$
Condition (i) implies that $\beta_1+N-(N-2)(p_1+q_1)>0$ and condition (ii) implies that $\beta_2+N-(N-2)(p_2+q_2)>0$.
Hence, in each of these cases, (SH)-(b) holds, and so (\ref{sys 1.1})  has no positive supersolution. Thus, the Proposition holds in case (i) and case (ii).

\vskip6pt
When $p_1+q_1=p_{\beta_1}^*$ and  $p_2+q_2=p_{\beta_2}^*$, we see that
 $$
 h_1(x,t|x|^{2-N},t|x|^{2-N})|x|^N=t^{p_1+q_1} (\ln |x|)^{\tau_1},\quad \forall\, |x|>e,
 $$
 $$
 h_2(x,t|x|^{2-N},t|x|^{2-N})|x|^N=t^{p_2+q_2} (\ln |x|)^{\tau_2},\quad\forall\, |x|>e.
 $$

Let $\sigma_3=-\tau_1$  if $-1<\tau_1 <0$ and $\sigma_4=-\tau_2$ if $-1<\tau_2<0$;
otherwise, we let $\sigma_3=\sigma_4=\frac12$. Then for $|x|>e$, we find that
 $$
 h_1  (x,\,t|x|^{2-N},\, t|x|^{2-N})  |x|^{N}  (\ln|x|)^{\sigma_3}\ge t^{p_1+q_1},\ \
 h_2  (x,\,t|x|^{2-N},\, t|x|^{2-N})  |x|^{N}  (\ln|x|)^{\sigma_4}\ge t^{p_2+q_2}.
 $$

Furthermore, if $p_1>1$, for any $t>0$ we obtain        that
 \begin{eqnarray*}
 \frac{h_1 (x,\,t|x|^{2-N}(\ln|x|)^{\sigma_5},\,t|x|^{2-N})}{ t|x|^{-N}  (\ln|x|)^{\sigma_5}} = t^{p_1+q_1-1}(\ln|x|)^{\tau_1 +\sigma_5(p_1-1)}\to +\infty, \ \ {\rm as}\ |x|\to+\infty
\end{eqnarray*}
by choosing $\sigma_5 >0$ so that $\tau_1 +\sigma_5(p_1-1) >0$.

If $q_2>1$, for any $t>0$ we  derive that
 \begin{eqnarray*}
 \frac{h_2 (x,\,t|x|^{2-N} ,\,t|x|^{2-N}(\ln|x|)^{\sigma_6})}{ t|x|^{-N}  (\ln|x|)^{\sigma_6}} =s^{p_2+q_2-1}(\ln|x|)^{\tau_2 + \sigma_6(q_2-1)}\to +\infty, \ \ {\rm as}\ |x|\to+\infty
 \end{eqnarray*}
by choosing $\sigma_6>0$ so that $\tau_2+\sigma_6(q_2-1) >0$.

If $p_1=q_2=1$, then for any $t>0$, it holds   that
 $$
 \frac{h_1 (x,\,t|x|^{2-N}(\ln|x|)^{\sigma_5},\,t|x|^{2-N}(\ln|x|)^{\sigma_6})}{ t |x|^{-N}  (\ln|x|)^{\sigma_5}}
 =t^{q_1}(\ln|x|)^{\tau_1 +\sigma_6q_1},\quad \forall\, |x|>e,
 $$
 $$
 \frac{h_2 (x,\,t|x|^{2-N}(\ln|x|)^{\sigma_5},\,t|x|^{2-N}(\ln|x|)^{\sigma_6})}{t |x|^{-N}  (\ln|x|)^{\sigma_6}}
 =t^{q_1}(\ln|x|)^{\tau_2 +\sigma_5p_2},\quad \forall\, |x|>e.
 $$
Since $p_2,\,q_1\geq0$ and $p_2+q_1>0$, we may  choose either $\sigma_5>0 $ or $ \sigma_6>0$  such that
 $$
 \tau_2 +\sigma_5p_2>0\ \ \mbox{or}\ \ \tau_1 +\sigma_6q_1>0.
 $$

According to the above analysis, it is easily seen that (SH)-(b1)(b2)(b3) are fulfilled
when one of the cases (iii), (iv) and (v) holds. Thus, Theorem \ref{teo 3} applies to assert that
(\ref{sys 1.1})  has no positive supersolution in each of such cases.
The proof is complete. \hfill$\Box$

%\begin{remark}
%We note that  the nonexistence of positive bounded supersolution to  Problem (\ref{sys 1.1}) holds true under assumption Proposition \ref{cr 2-3}
%without the restrictions  $p_1,\, q_2\ge1$, $p_2,\, q_1\ge0$, $\beta_1,\, \beta_2>-2$.

%\end{remark}

\vskip8pt

In what follows, we will establish the existence of positive supersolutions to system (\ref{sys 1.1}).
For sake of convenience, denote
 \begin{equation}\label{s-3-2}
 \sigma_1=\frac{(\tau_2+1)q_1-(\tau_1+1)(q_2-1)}{(p_1-1)(q_2-1)-p_2q_1} \quad {\rm and}\quad
  \sigma_2=\frac{(\tau_1+1) p_2 -(\tau_2+1)(p_1-1)}{(p_1-1)(q_2-1)-p_2q_1}.
\end{equation}

\begin{proposition}\label{cr 2-6-} Assume that
 $$
 h_1(x,u,v)=|x|^{\beta_1} (\ln |x|)^{\tau_1}u^{p_1}v^{q_1}, \ \ h_2(x,u,v)=|x|^{\beta_2} (\ln |x|)^{\tau_2}u^{p_2}v^{q_2}.
 $$
Problem (\ref{sys 1.1})  has a positive supersolution if $\beta_1,\, \beta_2\in\mathbb{R}$, and
 $$
 p_1+q_1\geq p_{\beta_1}^*,\ \ p_2+q_2\geq p_{\beta_2}^*,\ \ \sigma_1>0,\ \ \sigma_2>0.
 $$

\end{proposition}

\noindent{\bf Proof.} Set
 $$
 w_i(x)=|x|^{2-N} (\ln|x|)^{\sigma_i}\ \ \ \mbox{with}\ \sigma_i>0,\ i=1,2.
 $$
If $|x|>\max\big\{1,\  e^{\frac{2(\max\{\sigma_1,\sigma_2\}-1)}{N-2}}\big\}$,
it follows from (\ref{3.0}) that
 \begin{eqnarray*}
 -\Delta( tw_i(x)) &=& t\sigma_i(N-2)|x|^{-N} (\ln |x|)^{\sigma_i-1}-t\sigma_i(\sigma_i-1) |x|^{-N} (\ln |x|)^{\sigma_i-2} \nonumber \\
 &\ge &  \frac12t\sigma_i (N-2)|x|^{-N} (\ln |x|)^{\sigma_i-1},\ \ i=1,2
 \end{eqnarray*}
and since $p_1+q_1\geq p^*_{\beta_1}$,
 \begin{eqnarray*}
 h_1(x,tw_1,tw_2) &=& t^{p_1+q_1} |x|^{\beta_1} (\ln |x|)^{\tau_1}  |x|^{(2-N)(p_1+q_1)} (\ln|x|)^{\sigma_1p_1+\sigma_2q_1}   \\
 &\geq&   t^{p_1+q_1} |x|^{-N}  (\ln|x|)^{\tau_1+\sigma_1p_1+\sigma_2q_1}.
 \end{eqnarray*}
Similarly, by $p_2+q_2\geq p^*_{\beta_2}$, we have that
 \begin{eqnarray*}
 h_2(x,tw_1,tw_2)  \geq t^{p_2+q_2} |x|^{-N}  (\ln|x|)^{\tau_2+\sigma_1p_2+\sigma_2q_2}.
 \end{eqnarray*}
In view of the definitions of $\sigma_1,\sigma_2$,
 \begin{equation}\label{sigma}
 \tau_1+\sigma_1p_1+\sigma_2q_1=\sigma_1-1,\quad \tau_2+\sigma_1p_2+\sigma_2q_2=\sigma_2-1.
 \end{equation}
By choosing
 $$
 t=\min\left\{\left[\frac12\sigma_1 (N-2)\right]^{\frac1{p_1+q_1}},\ \ \left[\frac12\sigma_2 (N-2)\right]^{\frac1{p_2+q_2}}\right\},
 $$
we then see that $(tw_1,\, tw_2)$ is a supersolution of (\ref{sys 1.1}).\hfill$\Box$

\begin{remark}\label{re-b} Concerning the condition $\sigma_1,\,\sigma_2>0$ imposed
in Proposition \ref{s-3-2}, we want to make some comments as follows.

 \begin{itemize}
 \item[{\rm (i)}] When $p_1=q_2=1$, then
 $$
 \sigma_1=-\frac{\tau_2+1}{p_2}>0\quad {\rm and}\quad\sigma_2=-\frac{\tau_1+1}{q_1}>0
 $$
holds if $\tau_1,\,\tau_2<-1,\,p_2,\,q_1>0$.

 \item[{\rm (ii)}] When  $q_2=1$ and $\tau_1=-1$, we have
 $$
 \sigma_1=-\frac{\tau_2+1}{p_2}>0\quad {\rm and}\quad\sigma_2=\frac{(\tau_2+1)(p_1-1)}{p_2 q_1}>0
 $$
if $\tau_2<-1,\,p_1<1,\,p_2,\,q_1>0.$

 \end{itemize}

\end{remark}

%%%%%%%%%%%%%%%%%%%%%%%%%%%%%%%%%%%%%%%%%%%%%%%%%%%%%%%%%%%%%%%%%%%%%%

%%%%%%%%%%%%%%%%%%%%%%%%%%%%%%%%%%%%%%%%%%%%%%%%%%%%%%%
%%% Acknowledgements. 致谢
%%%%%%%%%%%%%%%%%%%%%%%%%%%%%%%%%%%%%%%%%%%%%%%%%%%%%%%
\Acknowledgements{H. Chen  is supported by NSF of China (No. 11726614, 11661045),
and Jiangxi Provincial Natural Science Foundation (No. 20161ACB20007).
R. Peng is supported by NSF of China (No. 11671175, 11571200),
the Priority Academic Program Development of Jiangsu Higher Education Institutions,
Top-notch Academic Programs Project of Jiangsu Higher Education Institutions (No. PPZY2015A013)
and Qing Lan Project of Jiangsu Province.
F. Zhou is supported by NSF of China (No. 11726613, 11271133, 11431005)
and STCSM (No. 13dZ2260400). The authors would like to thank Professor Yehuda Pinchover for bringing to our attention
some existing work on the linear Hardy potential problem (\ref{eq 1.1-ab}) such as \cite{DFP,P,P1} and
the two referees for their careful reading and valuable suggestions, which helped to
improve the exposition of the paper.
}

%%%%%%%%%%%%%%%%%%%%%%%%%%%%%%%%%%%%%%%%%%%%%%%%%%%%%%%
%%% Conflict of interest. 作者利益声明
%%%%%%%%%%%%%%%%%%%%%%%%%%%%%%%%%%%%%%%%%%%%%%%%%%%%%%%
%\InterestConflict

%%%%%%%%%%%%%%%%%%%%%%%%%%%%%%%%%%%%%%%%%%%%%%%%%%%%%%%
%%% Supplements. 补充材料, 非必选
%%%%%%%%%%%%%%%%%%%%%%%%%%%%%%%%%%%%%%%%%%%%%%%%%%%%%%%
%\Supplements{}

%%%%%%%%%%%%%%%%%%%%%%%%%%%%%%%%%%%%%%%%%%%%%%%%%%%%%%%
%%% Reference section. 参考文献
%%% citation in the content using "some words~\cite{1,2}".
%%% ~ is needed to make the reference number is on the same line with the word before it.
%%%%%%%%%%%%%%%%%%%%%%%%%%%%%%%%%%%%%%%%%%%%%%%%%%%%%%%

%%%%%%%%%%%%%%%%%%%%%%%%%%%%%%%%%%%%%%%%%%%%%%%%%%%%%%%
%%% Appendix sections. 附录章节, 非必选
%%%%%%%%%%%%%%%%%%%%%%%%%%%%%%%%%%%%%%%%%%%%%%%%%%%%%%%

\end{document}